\documentclass[11pt]{amsart}

\usepackage{dsfont,bbm,color,bm,natbib,multirow, graphicx}
\usepackage{amssymb,amsmath,amsthm}
\numberwithin{equation}{section}
\newtheorem{theorem}{Theorem}[section]
\newtheorem{lemma}[theorem]{Lemma}
\newtheorem{remark}[theorem]{Remark}

\begin{document}

\newcommand{\us}{\underline{s}}
\newcommand{\um}{\underline{m}}
\newcommand{\E}{\mathbb{E}}
\newcommand{\var}{\mathbb{V}ar}
\newcommand{\cov}{\mathbb{C}ov}
\renewcommand{\P}{\mathbb{P}}
\renewcommand{\a}{\alpha}
\renewcommand{\b}{\beta}
\newcommand{\Be}{\pmb\beta_0}
\renewcommand{\k}{\kappa}
\renewcommand{\t}{\tau}
\newcommand{\ep}{\varepsilon}
\newcommand{\tep}{\tilde\varepsilon}
\newcommand{\hep}{\hat\varepsilon}
\newcommand{\noin}{\noindent}
\newcommand{\De}{\Delta}
\newcommand{\la}{\lambda}
\newcommand{\La}{\Lambda}
\newcommand{\de}{\delta}
\renewcommand{\th}{\theta}
\newcommand{\si}{\sigma}
\newcommand{\Ga}{\Gamma}
\newcommand{\ga}{\gamma}
\newcommand{\pmu}{\pmb\mu}
\newcommand{\Si}{\Sigma}

\renewcommand{\o}{\omega}
\renewcommand{\O}{\Omega}
\newcommand{\ol}{\overline}
\newcommand{\bbI}{{\mathbb I}}
\newcommand{\bbC}{{\mathbb C}}
\newcommand{\bbR}{{\mathbb R}}
\newcommand{\cC}{{\mathcal C}}
\newcommand{\cD}{{\mathcal D}}
\newcommand{\cG}{{\mathcal G}}
\newcommand{\cI}{{\mathcal I}}
\newcommand{\cJ}{{\mathcal J}}
\newcommand{\cK}{{\mathcal K}}
\newcommand{\cQ}{{\mathcal Q}}

\newcommand{\frB}{\mathfrak{B}}
\newcommand{\frm}{\mathfrak{m}}
\newcommand{\frt}{\mathfrak{t}}
\newcommand{\cL}{{\mathcal L}}
\newcommand{\cH}{{\mathcal H}}
\newcommand{\cR}{{\mathcal R}}
\newcommand{\cT}{{\mathcal T}}
\newcommand{\cU}{{\mathcal U}}
\newcommand{\cV}{{\mathcal V}}
\newcommand{\vp}{\varpi}
\newcommand{\trvp}{\varpi^{tr}}
\newcommand{\evp}{\varpi^{\E}}

\newcommand{\A}{{\bf A}}
\newcommand{\C}{{\bf C}}
\newcommand{\B}{{\bf B}}
\newcommand{\D}{{\bf D}}
\newcommand{\e}{{\bf e}}
\newcommand{\F}{{\bf F}}
\newcommand{\G}{{\bf G}}
\renewcommand{\H}{{\bf H}}
\newcommand{\I}{{\bf I}}
\newcommand{\J}{{\bf J}}
\newcommand{\K}{{\bf K}}
\renewcommand{\L}{{\bf L}}
\newcommand{\M}{{\bf M}}
\newcommand{\N}{{\bf N}}
\newcommand{\Q}{{\bf Q}}
\renewcommand{\O}{{\bf O}}
\renewcommand{\r}{{\bf r}}
\renewcommand{\S}{{\bf S}}
\newcommand{\T}{{\bf T}}
\newcommand{\U}{{\bf U}}
\newcommand{\V}{{\bf V}}
\newcommand{\W}{{\bf W}}
\newcommand{\X}{{\bf X}}
\newcommand{\x}{{\bf x}}
\newcommand{\Y}{{\bf Y}}
\newcommand{\Z}{{\bf Z}}
\newcommand{\bz}{{\bf z}}
\renewcommand{\(}{\left(}
\renewcommand{\)}{\right)}
\newcommand{\lj}{\left|}
\newcommand{\rj}{\right|}
\newcommand{\lb}{\label}
\newcommand{\no}{\nonumber\\}

\newcommand{\hi}{H\"{o}lder's inequality}


\title[ On testing the equality of mean vectors]{On testing the equality of high dimensional mean vectors with unequal covariance matrices}

\author{Jiang Hu, \ \ Zhidong Bai, \ \ Chen Wang, \ \ Wei Wang
}
\thanks{J. Hu was partially supported by CNSF 11301063. Z. D. Bai was partially supported by CNSF 11171057. }

\address{KLASMOE and School of Mathematics \& Statistics, Northeast Normal University, Changchun, P.R.C., 130024.}
\email{huj156@nenu.edu.cn}

\address{KLASMOE and School of Mathematics \& Statistics, Northeast Normal University, Changchun, P.R.C., 130024.}
\email{baizd@nenu.edu.cn}

\address{Department of Statistics and Applied Probability, National University of Singapore, Singapore, 117546.}
\email{stawc@nus.edu.sg}
\address{KLASMOE and School of Mathematics \& Statistics, Northeast Normal University, Changchun, P.R.C., 130024.}
\email{wangw044@nenu.edu.cn}
%
%
%


\maketitle

\begin{abstract}
In this article, we focus on  the problem of  testing the equality of several high dimensional mean vectors with unequal covariance matrices. This is one of the most  important problems
in multivariate statistical analysis and there have been various tests proposed in the literature. Motivated by \citet{BaiS96E} and \cite{ChenQ10T},
we introduce a test statistic and derive the asymptotic distributions under the null and the alternative hypothesis.
 In addition, it is compared with a test statistic recently proposed by \cite{SrivastavaK13Ta}.  It is shown that our test statistic performs much better especially in the large dimensional case.
\end{abstract}

\section{ Introduction.}

In the last three decades, more and more large dimensional data sets appear in scientific
research. When the dimension of data or the number of parameters becomes large, the classical methods could reduce statistical efficiency significantly. In order to analyze those large data sets, many new statistical techniques, such as large dimensional multivariate
statistical analysis  based on the random matrix theory, have been developed. In this article, we consider the problem of testing the equality of several  high dimensional mean vectors with unequal covariance matrices, which is also called multivariate analysis of variance (MANOVA) problem. This problem is one of the most common multivariate statistical procedures in the social science, medical science, pharmaceutical science and genetics.  For example, a kind of disease may have several treatments. In the past, doctors only concern which treatments can cure the disease, and the standard clinical cure is low dimension. However, nowadays researchers  want to know whether the treatments alter  some of the  proteins or genes, thus then the high dimensional MANOVA is needed.

 Suppose there are $k(k\ge 3)$ groups and $X_{i1},\dots,X_{in_i}$ are $p$-variate independent and identically distributed (i.i.d.) random samples vectors from the $i$-th group, which have mean vector $\mu_i$ and covariance matrix $\Si_i$. We consider the problem of testing the hypothesis:
\begin{align}\lb{h}
  H_0:\mu_1=\cdots=\mu_k \quad \mbox{vs}\quad H_1:\exists i\neq j,~~ \mu_i\neq\mu_j.
\end{align}
Notice that here we do not need normality assumption. 
  The MANOVA problem has been discussed intensively  in the literature about multivariate statistic analysis. For example, for normally distributed groups, when the total sample size $n=\sum_{i=1}^kn_i$ is considerably larger than the dimension $p$, statistics that have been commonly used are likelyhood ratio test statistic \citep{Wilks32C},  generalized $T^2$ statistic \citep{Lawley38G,Hotelling47M} and  Pillai statistic \citep{Pillai55S}.
When $p$ is larger than the sample size $n$,  \cite{Dempster58H,Dempster60S} firstly considered this problem  in the case of two sample problem.  Since then, more  high dimensional tests have been proposed  by \citet{BaiS96E,FujikoshiH04A,SrivastavaF06M,Srivastava07M, Schott07S,SrivastavaD08T,Srivastava09T, SrivastavaY10T,ChenQ10T,SrivastavaK11S, SrivastavaK13T}. And recently, \cite{CaiX14H} proposed a statistic to  test the equality of multiple high-dimensional mean vectors under common covariance matrix.
Also, one can refer to the book \citep{FujikoshiU11M} for more details.

The statistic of testing (\ref{h}) we proposed in this article is
motivated by \citet{BaiS96E} and \cite{ChenQ10T}.
Firstly, let us review the two test statistics briefly.
For  $k=2$ and $\Si_1=\Si_2=\Si$, \citet{BaiS96E} proposed the test statistic
\begin{align}\lb{bs}
  T_{bs}=(\bar{X}_1-\bar{X}_2)'(\bar{X}_1-\bar{X}_2)-\frac{n_1+n_2}{n_1n_2}trS_n,
\end{align}
and showed that  under some conditions, as  $\min\{p,n_1,n_2\}\to\infty$, $p/(n_1+n_2)\to y>0$ and $n_1/(n_1+n_2)\to\kappa \in(0,1)$
\begin{align*}
  \frac{T_{bs}-\|\mu_1-\mu_2\|^2}{\sqrt{Var(T_{bs})}}\stackrel{d}{\to} N(0,1).
\end{align*}
Here
\begin{align*}
    \bar X_i=\frac{1}{n_i}\sum_{j=1}^{n_i}X_{ij},\quad S_n=\frac{1}{n_1+n_2-2}\sum_{i=1}^2\sum_{j=1}^{n_i}(X_{ij}-\bar X_i)'(X_{ij}-\bar X_i)
\end{align*} and
\begin{align*}
Var(T_{bs})=\frac{2(n_1+n_2)^2(n_1+n_2-1)}{n_1^2n_2^2(n_1+n_2-2)}tr\Si^2(1+o(1)).
\end{align*}
 In addition, Bai and Saranadasa gave a ratio-consistent estimator of $tr\Si^2$ (􏰻in the sense
that $\widehat{tr\Si^2}/tr\Si^2\to1$),  􏱯 􏰿􏰼􏱆that was
\begin{align*}
\widehat{tr\Si^2}=\frac{(n_1+n_2-2)^2}{(n_1+n_2)(n_1+n_2-3)}\(trS_n^2-\frac{1}{n_1+n_2-2}(trS_n)^2\).
\end{align*}

If $\Si_1\neq\Si_2$,    \cite{ChenQ10T}  gave a test statistic
\begin{align*}
  T_{cq}=\frac{\sum_{i\neq j}^{n_1}X_{1i}'X_{1j}}{n_1(n_1-1)}+\frac{\sum_{i\neq j}^{n_2}X_{1i}'X_{2j}}{n_2(n_2-1)}-2\frac{\sum_{i=1}^{n_1}\sum_{j=1}^{n_2}X_{1i}'X_{2j}}{n_1n_2},
\end{align*}
which can be expressed as
\begin{align}\lb{cq}
  T_{cq}=(\bar{X}_1-\bar{X}_2)'(\bar{X}_1-\bar{X}_2)-n_1^{-1}trS_1-n_2^{-1}trS_2.
\end{align}
Here and throughout this paper, the sample covariance matrix of the $i$-th group is denoted as
\begin{align*}
S_i=\frac{1}{n_i-1}\sum_{j=1}^{n_i}(X_{ij}-\bar X_i)'(X_{ij}-\bar X_i).
\end{align*}
Also  they proved that under some conditions
\begin{align*}
  \frac{T_{cq}-\|\mu_1-\mu_2\|^2}{\sqrt{Var(T_{cq})}}\stackrel{d}{\to} N(0,1)
\end{align*}
where
\begin{align*}
  Var(T_{cq})=\(\frac{2}{n_1(n_1-1)}tr(\Si_1^2)+\frac{2}{n_2(n_2-1)}tr(\Si_2^2)+\frac{4}{n_1n_2}tr(\Si_1\Si_2)\)(1+o(1)).
\end{align*}
And then  \cite{ChenQ10T} gave the  ratio-consistent estimators of $tr\Si_i^2$ and $tr(\Si_1\Si_2)$, that were
\begin{align}\lb{est2}
  \widehat{tr(\Si_i^2)}=\frac{1}{n_i(n_i-1)}tr\(\sum_{j\neq k}^{n_i}(X_{ij}-\bar X_{i(j,k)})X_{ij}'(X_{ik}-\bar X_{i(j,k)})X_{ik}'\)
\end{align}
and
\begin{align}\lb{est3}
  \widehat{tr(\Si_1\Si_2)}=\frac{1}{n_1n_2}tr\(\sum_{l=1}^{n_1}\sum_{k=1}^{n_2}(X_{1l}-\bar X_{1(l)})X_{1l}'(X_{2k}-\bar X_{2(k)})X_{2k}'\).
\end{align}
Here $\bar X_{i(j,k)}$ is the $i$-th sample mean after excluding $X_{ij}$ and $X_{ik}$, and $\bar  X_{i(l)}$ is
the $i$-th sample mean without $X_{il}$.

When $\Si_1=\Si_2$, it is apparent that the test statistic proposed by \cite{ChenQ10T} reduces to the one obtained by  \cite{BaiS96E}. Compared to \cite{BaiS96E}, \cite{ChenQ10T} generalized the test to the case when $\Si_1\neq \Si_2$, and  used   different estimators of the variance. This is indeed a significant improvement to remove the assumption $\Si_1= \Si_2$, because such an assumption is hard to verify for high-dimensional data.
Thus based on these properties,  we propose a statistic of  testing the equality of more than two high dimensional mean vectors with unequal covariance matrices.

We assume the following general multivariate model:
\begin{itemize}
  \item[(a):] $
  X_{ij}=\Ga_iZ_{ij}+\mu_i,\quad \mbox{for } i=1,\dots k, j=1\dots,n_i,
$
where $\Ga_i$ is a $p\times m$ matrix for some $m\geq p$ such that $\Ga_i\Ga_i'=\Si_i$, and $\{Z_{ij}\}_{j=1}^{n_i}$ are $m$-variate i.i.d. random vectors satisfying $E(Z_{ij})=0$ and $Var(Z_{ij})=I_m$, the $m\times m$ identity matrix;
  \item[(b):] $Z_{ij}=(z_{ij1},\dots,z_{ijm})'$, with
$
  E(z_{ijl_1}^{\a_1}z_{ijl_2}^{\a_2}\dots z_{ijl_q}^{\a_q})= E(z_{ijl_1}^{\a_1})E(z_{ijl_2}^{\a_2})\dots E( z_{ijl_q}^{\a_q})
$ and $E(z_{ijk}^4)<\infty$,
for a positive integer $q$ such that $\sum_{l=1}^q\a_l\leq 8$ and $l_1\neq l_2 \neq\dots\neq l_q$;
  \item[(c):]$
\frac{n_i}{n}\to k_i\in(0,1)\quad i=1,\dots k,  \mbox{~as~~} n\to\infty.
$ Here  $n=\sum_{i=1}^kn_i$;

    \item[(d):]$tr(\Sigma_l\Sigma_d\Sigma_l\Sigma_h
)=o[tr(\Sigma_l\Sigma_d)tr(\Sigma_l\Sigma_h
)]\qquad d,l,h\in\{1,2,\dots, k\}; $
\item[(e):]$
(\mu_d-\mu_l)'
\Sigma_d(\mu_d-\mu_h)=o[n^{-1}tr\{(\sum_{i=1}^k\Sigma_i)^2\}],\quad d,l,h\in\{1,2,\dots, k\}
$.
\end{itemize}
 It should be noted that all random variables and parameters here and later
depend on $n$. For simplicity we omit the subscript $n$ from all random variables
except those statistics defined later.

Now we construct our test. Consider the statistic
\begin{eqnarray*}
  T_n^{(k)}&=&\sum_{i<j}^k(\bar X_i-\bar X_j)'(\bar X_i-\bar X_j)-(k-1)\sum_{i=1}^kn_i^{-1}tr S_i\\
  &=&(k-1)\sum_{i=1}^k\frac1{n_i(n_i-1)}\sum_{k_1\ne k_2}X_{ik_1}'X_{ik_2}-\sum_{i<j}^k\frac2{n_in_j}\sum_{k_1,k_2}X_{ik_1}'X_{jk_2}.
\end{eqnarray*}
When $k=2$, apparently $T_n^{(2)}$ is the Chen-Qin test statistic. Next we will calculate the mean and variance of $T_n^{(k)}$. Unlike the method used in  \cite{ChenQ10T}, we give a much simpler procedure.
From $X_{ij}=\Ga_iZ_{ij}+\mu_i$,  we can rewrite $T_n^{(k)}-\sum_{i<j}\|\mu_i-\mu_j\|^2$ as $T_1^{(k)}+T_2^{(k)}$, where
\begin{eqnarray*}
T_1^{(k)}&=&(k-1)\sum_{i=1}^k\frac{1}{n_i(n_i-1)}\sum_{k_1\ne k_2}Z_{ik_1}'\Ga_i'\Ga_i Z_{ik_2}-\sum_{i<j}^k\frac2{n_in_j}\sum_{k_1,k_2}Z_{ik_1}'\Ga_i'\Ga_jZ_{jk_2}\\
T_2^{(k)} &=&\sum_{i=1}^k\frac2{n_i}(k\mu_i-\sum_{j=1}^k\mu_j)'\sum_{k}\Ga_iZ_{ik_1}.
\end{eqnarray*}
Thus we can show immediately  that
\begin{align*}
  E(T_n^{(k)})=\sum_{i<j}^k\| \mu_i- \mu_j\|^2
\end{align*}
and
\begin{eqnarray*}
  Var(T_n^{(k)})&=&\sum_{i=1}^k\frac{2(k-1)^2}{n_i(n_i-1)}tr(\Sigma_i^2)
+\sum_{i< j}^k\frac{4}{n_in_j}tr(\Sigma_i\Sigma_j)\\
&&+4\sum_{i=1}^k\frac1{n_i}\left(\sum_{j=1}^k\mu_j-k\mu_i\right)'\Sigma_i\left(\sum_{j=1}^k\mu_j-k\mu_i\right).
\end{eqnarray*}
Then we have the following theorem:
\begin{theorem}\lb{th1}
  Under the assumptions (a)-(e), we obtain that as $p\to\infty$ and $n\to\infty$,
  \begin{align}\lb{lim1}
    \frac{T_n^{(k)}-\sum_{i<j}^k\| \mu_i- \mu_j\|^2}{\sqrt{Var(T_n^{(k)})}}\stackrel{d}{\to}N(0,1).
  \end{align}
\end{theorem}
 It is worth  noting that under $H_0$, assumption (e) is trivially satisfied and $E(T_n^{(k)})=0$. What is more, under $H_1$ and assumptions (a)-(e), $Var(T_n^{(k)})=(\si_n^{(k)})^2(1+o(1))$, where
 \begin{align*}
   (\si_n^{(k)})^2=\sum_{i=1}^k\frac{2(k-1)^2}{n_i(n_i-1)}tr(\Sigma_i^2)
+\sum_{i< j}^k\frac{4}{n_in_j}tr(\Sigma_i\Sigma_j).
\end{align*}
Then Theorem \ref{th1} is still true if the denominator of (\ref{lim1}) is replaced by $\si_n^{(k)}$. Therefore, to
complete the construction of our test statistic, we only need to find a ratio-consistent
estimator of $(\si_n^{(k)})^2$ and  substitute it into the denominator of (\ref{lim1}).
There are many estimators for $(\si_n^{(k)})^2$, and in this paper we choose two of them:
\begin{lemma}[UMVUE]\lb{th2}
 Under the assumptions (a)-(d), we obtain that as $p\to\infty$ and $n\to\infty$,
 \begin{align*}
 \frac{ \widehat{tr\Si_i^2}}{ tr\Si_i^2}\stackrel{p}{\to}1\mbox{\quad and\quad} \frac{ \widehat{tr(\Si_i\Si_j)}}{ tr(\Si_i\Si_j)}\stackrel{p}{\to}1
 \end{align*}
%
 where $i\neq j\in\{1,2,\dots, k\},$
\begin{align}\lb{est1}
 \widehat{tr\Si_i^2}=\frac{(n_i-1)^2}{(n_i+1)(n_i-2)}\(trS_i^2-\frac{1}{n_i-1}tr^2 S_i\)\end{align}
and\begin{align}\lb{est2}
\widehat{tr(\Si_i\Si_j)}=trS_iS_j.
\end{align} 

\end{lemma}

\begin{remark}
Under the normality assumption \eqref{est1} and \eqref{est2}
are uniformly minimum variance unbiased
estimators. The proof of this lemma was given in \cite{BaiS96E} and \cite{Srivastava09T}, and we omit it in this paper. 
\end{remark}

\begin{lemma}[Unbiased nonparametric estimators (UNE)]\label{leub}
 Under the assumptions (a)-(d), we obtain that as $p\to\infty$ and $n\to\infty$,
 \begin{align*}
 \frac{ \widehat{tr\Si_i^2}}{ tr\Si_i^2}\stackrel{p}{\to}1\mbox{\quad and\quad} \frac{ \widehat{tr(\Si_i\Si_j)}}{ tr(\Si_i\Si_j)}\stackrel{p}{\to}1
 \end{align*}
  where $i\neq j\in\{1,2,\dots, k\},$
\begin{align}\lb{esub1}
 \widehat{tr\Si_i^2}&=\frac{1}{(n_i)_6}\times\nonumber \\&\sum_{k_1,\dots,k_6\atop distinct }(X_{ik_1}-X_{ik_2})'(X_{ik_3}-X_{ik_4})(X_{ik_3}-X_{ik_5})'(X_{ik_1}-X_{ik_6})
 \end{align}
 and
 \begin{align}\lb{esub2}
 \widehat{tr(\Si_i\Si_j)}&=\frac1{(n_i)_3(n_j)_3}\times\nonumber\\&\sum_{k_1,k_2,k_3\ distict\atop
 k_4k_5,k_6\ distinct} (X_{ik_1}-X_{ik_2})'(X_{jk_4}-X_{jk_5})(X_{jk_4}-X_{jk_6})'(X_{ik_1}-X_{ik_3}).
\end{align}
Here $(n)_l=n(n-1)\cdots(n-l+1)$.
\end{lemma}

\begin{remark}
By Assumption $(a)$, the unbiasedness of estimators $\widehat{tr\Si_i^2}$ and $\widehat{tr\Si_i\Si_j}$ can be  easily proved and their  ratio-consistency can be found in \cite{LiC12T}.

\end{remark}

\begin{remark}
\cite{LiC12T} mentioned that the computation of  the estimators in Lemma \ref{leub} would be extremely heavy if the sample sizes  are very large. 
Thus to increase the computation speed, we simplify the estimator \eqref{esub1} further to:
  \begin{align*}\widehat{tr\Si_i^2}=&\frac{1}{n_i(n_i-3)}\|\Theta_i\|_2^2
-\frac2{n_i(n_i-2)(n_i-3)}\|\Theta_i\|_{1,2}^2
\\&+\frac1{n_i(n_i-1)(n_i-2)(n_i-3)}(\|\Theta_i\|_1)^2
\end{align*}
where $\Theta_i=X_i'X_i-Diag[X_i'X_i]$,  $X_i=(X_{i1},\dots X_{in_i})_{p\times n_i}$ and $Diag[X_i'X_i]$ is a diagonal matrix consisting of the diagonal elements of  $X_i'X_i$. Notice that for any matrix $A=(a_{ij})_{m\times n}$, the norm $\|\cdot\|_q$ is entrywise norm, i.e.,  $\|A\|_q=(\sum_{i=1}^m\sum_{j=1}^n|a_{ij}|^q)^{1/q}$  and the norm $\|\cdot\|_{p,q}$ is $L_{p,q}$ norm, i.e.,  $$\|A\|_{p,q}=(\sum_{i=1}^m(\sum_{j=1}^n|a_{ij}|^{p})^{q/p})^{1/q}.$$

What is more, from a direct calculation we can show that the estimator (\ref{esub2}) is exactly equal to the estimator \eqref{est2} in Lemma \ref{th2}.
That is because, 
$$(\ref{esub2}) =\frac{1}{(n_i-1)(n_j-1)}\sum_{k_1,k_4}(X_{ik_1}'X_{jk_4})^2 -\frac{1}{(n_i-1)n_j(n_j-1)}\sum_{k_1}^{n_i}(\sum_{k_4}^{n_j}X_{ik_1}'X_{jk_4})^2 $$
$$-\frac{1}{n_i(n_i-1)(n_j-1)}\sum_{k_4}^{n_j}(\sum_{k_1}^{n_i}X_{ik_1}'X_{jk_4})^2 +\frac{1}{n_i(n_i-1)n_j(n_j-1)}(\sum_{k_1, k_4}{X_{ik_1}}'X_{jk_4})^2 $$
$$=\frac{1}{(n_i-1)(n_j-1)}trX_iX_i'X_jX_j'-\frac{n_j}{(n_i-1)(n_j-1)}tr\bar{X}_i\bar{X}_i'{X}_j{X}_j' $$
$$-\frac{n_i}{(n_i-1)(n_j-1)} +\frac{n_in_j}{(n_i-1)(n_j-1)}tr\bar{X}_i\bar{X}_i'\bar{X}_j\bar{X}_j' =trS_iS_j.$$
Apparently,  using the simplified   formulas instead of  the original ones can make the   computation  much faster.

\end{remark}
Now,  by combining Theorem \ref{th1} and Lemma \ref{th2} (or Lemma \ref{leub}), we obtain our test statistic under $H_0$ and have the following theorem:
\begin{theorem}\lb{th3}
 Under $H_0$  and the assumptions (a)-(d), we obtain that as $p\to\infty$ and $n\to\infty$,
  \begin{align*}
   T_{our}={T_n^{(k)}}/{\hat\si_n^{(k)}}\stackrel{d}{\to}N(0,1),
  \end{align*}
  where $(\hat\si_n^{(k)})^2=\sum_{i=1}^k\frac{8}{n_i(n_i-1)}\widehat{ tr(\Sigma_i^2)}
+\sum_{i< j}^k\frac{4}{n_in_j}\widehat{ tr(\Sigma_i\Sigma_j)}$ with $\widehat{ tr(\Sigma_i^2)}$ and $\widehat{ tr(\Sigma_i\Sigma_j)}$ given in Lemma \ref{th2} or Lemma \ref{leub}.
\end{theorem}
\begin{remark}
When  the number of groups $k$ is small, the hypothesis $H_0$ can be considered as a multiple hypothesis of testing each two sample.  And the test for each sub-hypothesis can be tested by Chen and Qin (2010). However, for each sub-hypothesis, there is a test statistic of Chen and Qin (2010). The problem is how do we set up the critical value for the simultaneous test of the compound hypothesis $H_0$. In the literature, there is a famous Bonferroni correction method can be used. But it is well known that Bonfferoni correction is much conservative. Form this theorem, we can see that using our test, one may set up an asymptotically exact test.

\end{remark}

Due to Theorem \ref{th3}, the test with an $\a$ level of significance  rejects $H_0$ if $T_{our}>\xi_\a$ where $\xi_\a$ is the upper $\a$ quantile of $N(0,1)$. Next we will  discuss the power properties of the proposed test. Denote $\|\mu\|=\sum_{i<j}^k\| \mu_i- \mu_j\|^2$. From the above conclusions,  we can  easily obtain that
\begin{align*}
  T_{our}-\frac{\|\mu\|}{\sqrt{Var(T_n^{(k)})}}\stackrel{d}{\to}N(0,1).
\end{align*}
This implies
\begin{align*}
  \b_{nT}(\|\mu\|)=P_{H_1}(T_{our}>\xi_\a)=\Phi(-\xi_\a+\frac{\|\mu\|}{\si_n^{(k)}})+o(1),
\end{align*}
where $\Phi$ is the standard normal distribution function.

\section{Other tests  and simulations}
Due to the fact that the commonly used likelihood ratio test performs badly when dimension is large has been considered in a lot of literature such as \cite{BaiS96E,BaiJ09C,JiangJ12L,JiangY13C}, the discussion of the likelihood ratio test is left out  in this paper. Recently,
\cite{SrivastavaK13Ta} proposed a test statistic  of testing the equality of mean vectors of  several groups with a common unknown non-singular covariance matrix.  Denote ${\bf1}_r=(1,\dots,1)'$ as an $r$-vector with all the elements equal to one and define
$
  Y=(X_{11},\dots,X_{1n_1},\dots,X_{k1},\dots,X_{kn_k})
$, $L=(I_{k-1},-{\bf1}_{k-1})_{(k-1)\times k}$ and $$E=\left(
         \begin{array}{ccc}
           {\bf1}_{n_1} & \bf0 & \bf0 \\
            \bf0 & {\bf1}_{n_2} & \bf0  \\
           \vdots & \vdots & \vdots \\
           \bf0 & \bf0 & {\bf1}_{n_k} \\
         \end{array}
       \right)_{n\times k}.
$$
Then it is proposed that
\begin{align*}
  T_{sk}=\frac{tr(BD_S^{-1}-(n-k)p(k-1)(n-k-2)^{-1})}{\sqrt{2c_{p,n}(k-1)(trR^2-(n-k)^{-1}p^2)}},
\end{align*}
where $B=Y'E(E'E)^{-1}L'[L(E'E)^{-1}L']^{-1}L(E'E)^{-1}E'Y$, $D_S=Diag[(n-k)^{-1}Y(I_n-E(E'E)^{-1}E')Y]$, $R=D_S^{-1/2}Y(I_n-E(E'E)^{-1}E')YD_S^{-1/2}$ and $c_{p,n}=1+tr(R^2)/p^{3/2}.$ Notice that $Diag[A]$ denotes a diagonal matrix with the same diagonal elements as the diagonal elements of matrix $A$. Under the null hypothesis and the condition  $n=O(p^\de)$ with $\de>1/2$, $T_{sk}$ is asymptotically distributed as $N(0,1)$. That is as $n,p\to\infty$,
\begin{align*}
  P_{H_0}(T_{sk}>\xi_\a)\to\Phi(-\xi_\a).
\end{align*}

In this section we compare the performance of the proposed statistics $T_{our}$ and $T_{sk}$ in finite samples by simulation.   Notice that  the data is generated from the model
$$X_{ij}=\Ga_iZ_{ij}+\mu_i,\quad i=1,\dots,k,j=1,\dots,n_i.$$
where $\Ga_i$ is a $p\times p$  such that $\Ga_i^2=\Si_i$. 
 Here  $Z_{ij}=(z_{ij1},\dots,z_{ijp})'$ and $z_{ijk}$'s are independent random variables which are distributed as one of the following three
distributions:
\begin{align*}
  (\mathrm{i})~ N(0,1),\quad (\mathrm{ii})~(\chi^2_2-2)/2,\quad(\mathrm{iii})~(\chi^2_8-8)/4.
\end{align*}
For the covariance matrix $\Si_i,$ $i=1,2,3$, we consider the following two
cases:
\begin{itemize}
  \item[Case 1]: $\Si_i=\Ga_i=I_p$;
  \item[Case 2]: $\Si_i=\Ga_i^2=W_i\Psi_iW_i$, $W_i=Diag[w_{i1},\dots,w_{ip}]$, $w_{ij}=2*i+(p-j+1)/p$, $\Psi_i=(\phi^{(i)}_{jk})$, $\phi^{(i)}_{jj}=1$, $\phi^{(i)}_{jk}=(-1)^{j+k}(0.2\times i)^{|j-k|^{0.1}}$, $j\neq k$.
\end{itemize}

\begin{figure}[h]
\centering
\includegraphics[scale=0.33]{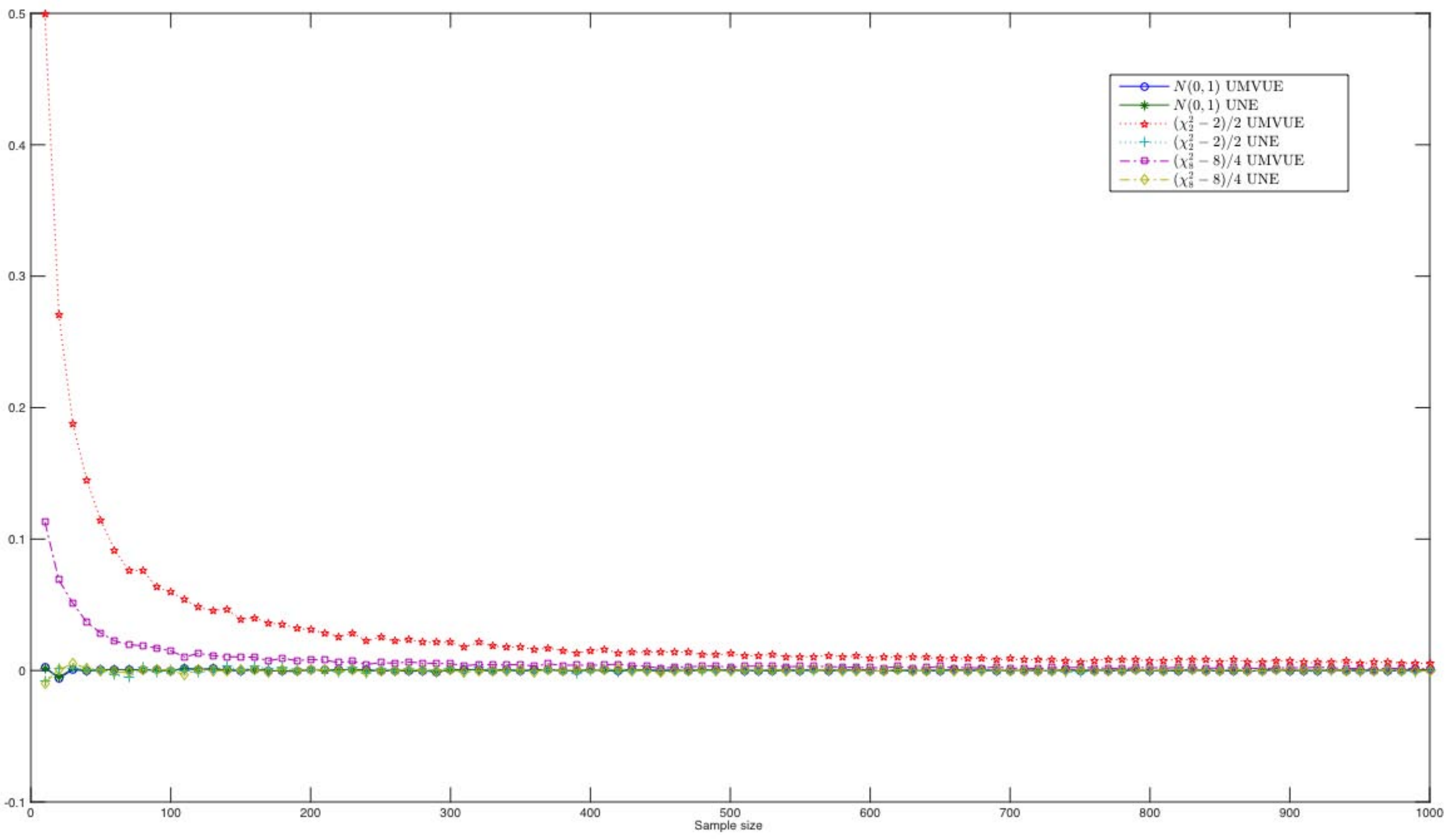}
\caption{Graph of the difference between the estimators $\frac{1}p\widehat{ tr\Sigma_1^2}$ and the true value  ${ \frac{1}p tr\Sigma_1^2}$ in Case 1, i.e., $\Si_1=I_p$. }\label{fig1}
\end{figure}

\begin{figure}[h]
\centering
\includegraphics[scale=0.33]{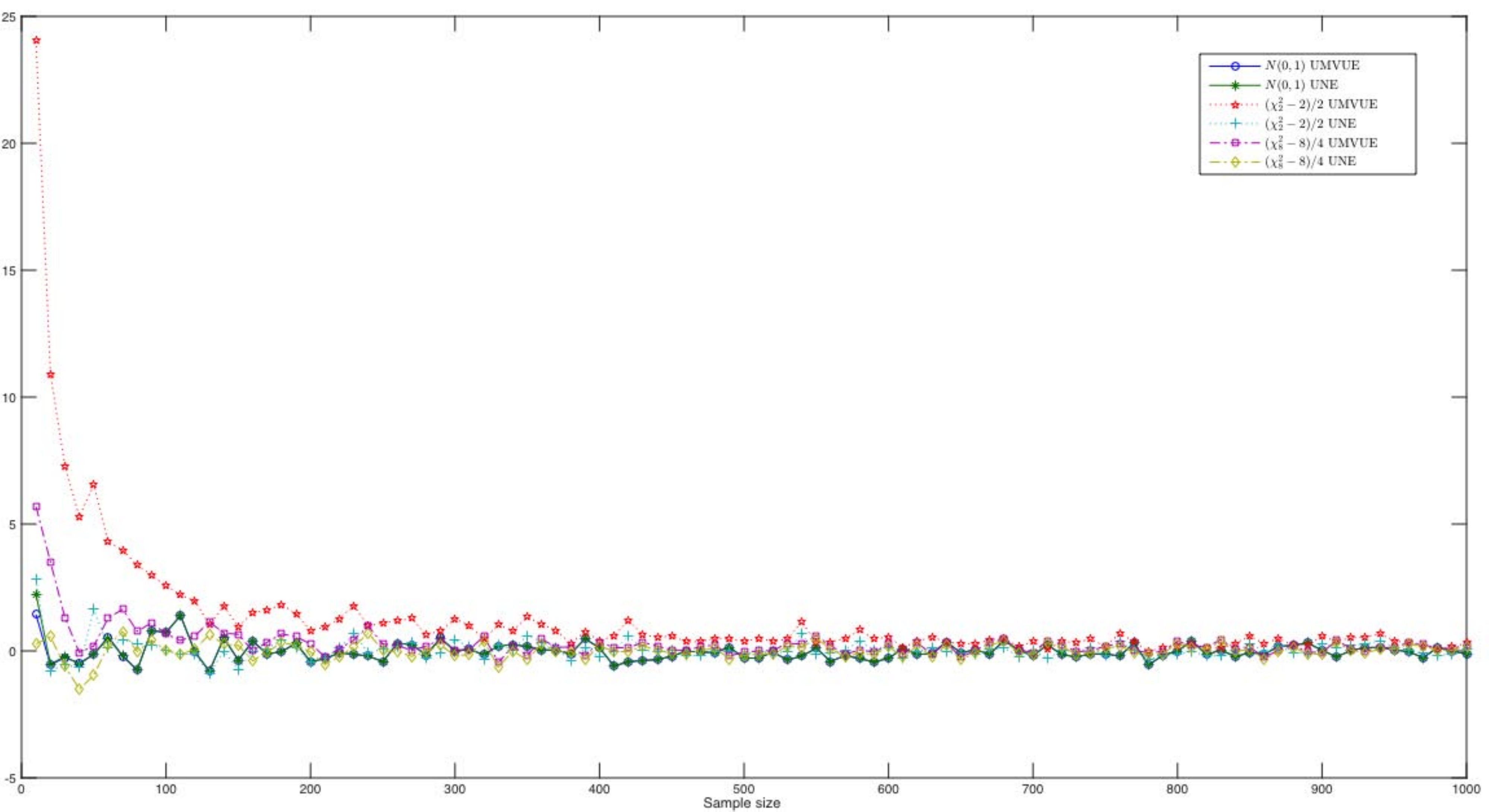}
\caption{Graph of the difference between the estimators $\frac{1}p\widehat{ tr\Sigma_1^2}$ and the true value  ${ \frac{1}p tr\Sigma_1^2}$ in Case 2, i.e., $\Si_1=W_1\Psi_1W_1$. }\label{fig2}
\end{figure}

We first compare the  convergence rates of the estimators  \eqref{est1} and \eqref{esub1} based on the above models, see Figure \ref{fig1} and Figure \ref{fig2}. Here the dimension $p=100$ and the sample sizes $n$ are from 10 to 1000. The results are based  on  1000 replications. From these two figures we can easily find that in both cases, the UNE \eqref{esub1} and UMVUE \eqref{est1} are almost the same if the data  sets come  from standard    normal distribution.  But UNE is much better  than UMVUE if the data sets come from $\chi^2$ distribution, especially when $n$ is small.  

\begin{figure}[h]
\centering
\includegraphics[scale=0.33]{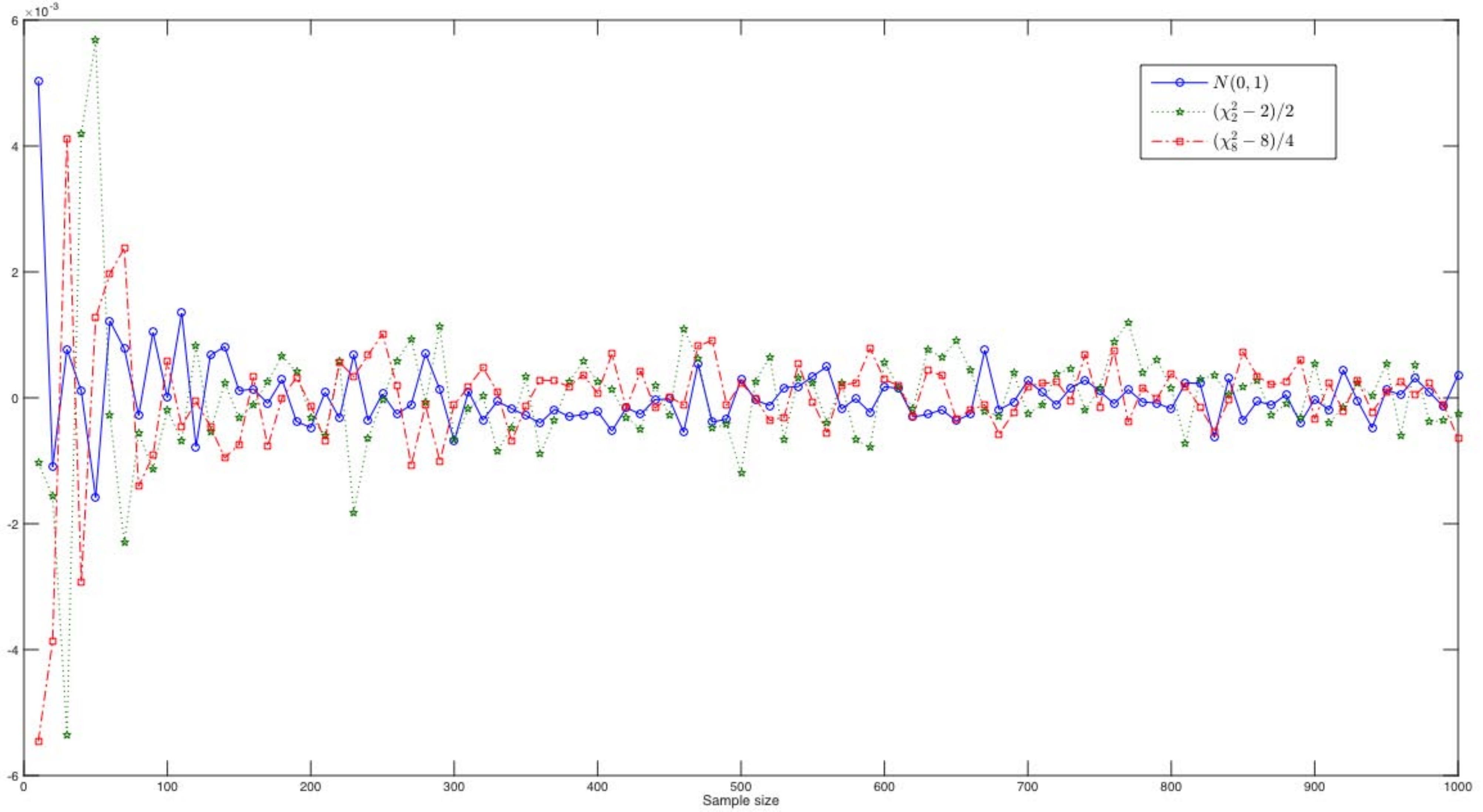}
\caption{Graph of the difference between the estimators $\frac{1}p\widehat{ tr\Sigma_1\Si_2}$ and the true value  ${ \frac{1}p tr\Sigma_1\Si_2}$ in Case 1, i.e., $\Si_1=\Si_2=I_p$. }\label{fig3}
\end{figure}

\begin{figure}[h]
\centering
\includegraphics[scale=0.33]{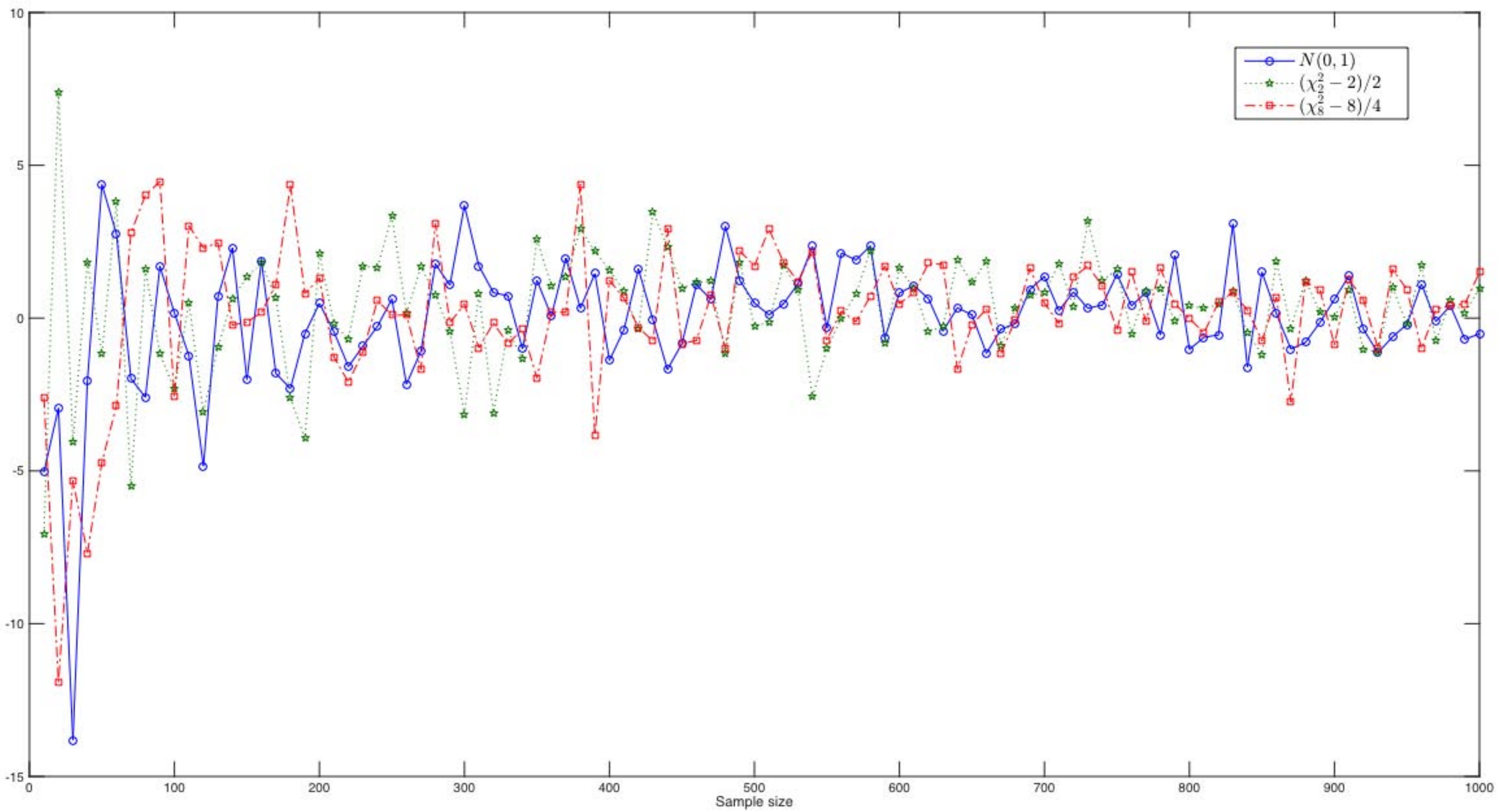}
\caption{Graph of the difference between the estimators $\frac{1}p\widehat{ tr\Sigma_1\Si_2}$ and the true value  ${ \frac{1}p tr\Sigma_1\Si_2}$ in Case 2, i.e., $\Si_1=W_1\Psi_1W_1$ and $\Si_2=W_2\Psi_2W_2$. }\label{fig4}
\end{figure}

Next let us see the performance  of the estimator $\widehat{tr\Sigma_i\Sigma_j}=trS_iS_j$ in Case 1 and Case 2 (see Figure \ref{fig3} and Figure  \ref{fig4}). Also the dimension $p=100$ and the sample sizes $n_1=n_2$ are from 10 to 1000. The results are based  on  1000 replications.  In Case 1, the estimator $\widehat{tr\Sigma_i\Sigma_j}=trS_iS_j$ performs very well  at all the three distributions. However in Case 2, we find that the convergence rate of this estimator is not very fast. Thus this will  be one of the  reasons that cause our test statistic $T_{our} $ to perform not well enough.


Now we examine the attained significance level (ASL) of the test statistics $T_{our}$ and $T_{sk}$ compared
to the nominal value $\a = 0. 05$, and then  examine their attained power.
The ASL is computed as $\hat\a=\#(T>\xi_{1-\a})/r$ where $T$ are values of the test statistic $T_{our}$ or $T_{sk}$ obtained from data simulated under $H_0$, $r$ is the number of replications and $\xi_{1-\a}$ is the $100 ( 1 -\a) \%$ quantile  of the standard normal distribution.
The attained power
of the test $T_{our}$ and $T_{sk}$ is also computed as
$\hat\b = \# ( T > \xi_{ 1-\a} )/ r$, where $ T$  are values of the test statistic $T_{our}$  or $T_{sk}$ computed from data simulated under the
alternative.

For simulation, we consider the problem of testing the equality of 3 mean vectors, that is, $k = 3$. Choose $p\in\{20,50,100,500,800\}$, $n_1=0.5\times n^*$, $n_2=n^*$, $n_3=1.5\times n^*$,  where $n^*\in\{20,50,100,200\}$. 
For the null hypothesis, without loss of generality we choose $\mu_1=\mu_2=\mu_3=\bf 0$.
For the alternative hypothesis, we choose $\mu_1=\bf0$, $\mu_2=(u_1,\dots,u_p)'$ and $\mu_3=-\mu_2$, where $u_i=(-1)^iv_i$ with $v_i$ are i.i.d. $U(0, a)$
which denotes uniform distribution with the support $( 0,a ) $.  Here in Case 1 we choose $a=0.1$ and in Case 2 we choose $a=0.2$, respectively.

The ASL and the powers are obtained based on 10,000
replications and  the $95\%$ quantile of the standard normal distribution is 1.64485. The  four tables report the ASL and the power in the {null} hypothesis and the alternative hypothesis of the two tests. For illustration,  in the tables we  respectively use the estimators proposed in Lemma \ref{th2} and Lemma \ref{leub} to obtain two different test statistics, $T_{our}^{umvue}$ and $T_{our}^{une}$. 
It is shown in Table 1 and Table 3 that the ASL of the proposed tests $T_{our}^{umvue}$ and $T_{our}^{une}$ approximate $\a = 0.05$ well in both cases,  and  $T_{our}^{une}$ is even better at nonnormal distributions.   Because it is shown that  UNE is   better than UMVUE if the data sets come from $\chi^2$ distribution, especially when $n$ is small.   But the ASL of test $T_{sk}$ in case 2 performs substantially worse. In addition, in Case 1 the  test $T_{sk}$ seems worse when dimension $p$ is much larger than the sample size $n^*$. This is probably because $T_{sk}$ is under common  covariance matrix assumption and needs condition $n=O(p^\de)$ with $\de>1/2$ to obtain the asymptotic distribution.
As reported
in  Table 2 and Table 4 , the powers of the  test $T_{our}^{une}$  perform better than $T_{sk}$ in Case 2 and worse in Case 1.  But actually in Case 1, when the dimension $p$ and  sample size  $n^*$ are large, the powers of the  test $T_{our}^{une}$ are also good enough. Thus when the dimension is much larger than the sample size, or the dimension and the sample size are both large, our test statistic is recommended, as it is more stable.

{
\begin{table}[p]
\caption{}
\begin{tabular}{ccccc|ccc|ccc}

 \multicolumn{11}{c}{ASL of $T_{our}$ and $T_{sk}$ in Case 1} \\
    \hline
\multirow{2}{*}{$p$} &\multirow{2}{*}{$n^*$}&\multicolumn{3}{c}{$z_{ijk}\stackrel{iid}{\sim}N(0,1)$}&\multicolumn{3}{c}{$z_{ijk}\stackrel{iid}{\sim}(\chi^2_2-2)/2$} &\multicolumn{3}{c}{$z_{ijk}\stackrel{iid}{\sim}(\chi^2_8-8)/4$}\\ \cline{3-11}
&&$T_{our}^{umvue}$&$T_{our}^{une}$&$T_{sk}$&$T_{our}^{umvue}$&$T_{our}^{une}$&$T_{sk}$&$T_{our}^{umvue}$&$T_{our}^{une}$&$T_{sk}$\\\hline
 20 & 20 & 0.0578 &0.0672 & 0.0505  & 0.0427 &0.0643& 0.0497 & 0.0549 &0.0655& 0.0531 \\
    & 50 & 0.0626 &0.0648 & 0.0480  & 0.0534  &0.0605& 0.0472 & 0.0578 &0.0594& 0.0489 \\
    & 100& 0.0630 &0.0626  & 0.0458 & 0.0543 &0.0614& 0.0466 & 0.0598 &0.0604& 0.0476 \\
    & 200& 0.0606 &0.0612 & 0.0434  & 0.0547 &0.0635& 0.0438 & 0.0587 &0.0586& 0.0441 \\
 50 & 20 & 0.0561 &0.0613& 0.0509  & 0.0383 &0.0604& 0.0471 & 0.0457 &0.0632& 0.0503 \\
    & 50 & 0.0551 &0.0573& 0.0454   & 0.0444 &0.0597& 0.0450 & 0.0538 &0.0572& 0.0443 \\
    & 100& 0.0596 &0.0608 & 0.0465  & 0.0528 &0.0545& 0.0490 & 0.0561 &0.0578& 0.0458 \\
    & 200& 0.0548 &0.0563& 0.0453   & 0.0528 &0.0580& 0.0453 & 0.0570 &0.0578& 0.0472 \\
100 & 20 & 0.0565 &0.0571& 0.0524  & 0.0346 &0.0572& 0.0428 & 0.0477 &0.0556& 0.0478 \\
    & 50 & 0.0551 &0.0573 & 0.0468   & 0.0471 &0.0599& 0.0476 & 0.0506 &0.0578& 0.0477 \\
    & 100& 0.0548 &0.0580& 0.0468   & 0.0499 &0.0570& 0.0471 & 0.0524 &0.0556& 0.0461 \\
    & 200& 0.0549 &0.0535& 0.0456   & 0.0541 &0.0599& 0.0455 & 0.0521 &0.0588& 0.0437 \\
500 & 20 & 0.0514 &0.0549& 0.0365  & 0.0315 &0.0543& 0.0351 & 0.0473 &0.0528& 0.0391 \\
    & 50 & 0.0535 &0.0553 & 0.0427  & 0.0422 &0.0515& 0.0434 & 0.0498 &0.0517& 0.0428 \\
    & 100& 0.0562 &0.0511& 0.0513   & 0.0457 &0.0510& 0.0423 & 0.0500 &0.0532& 0.0454 \\
    & 200& 0.0546 &0.0518& 0.0509   & 0.0497 &0.0477& 0.0473 & 0.0530 &0.0584& 0.0506 \\
800 & 20 & 0.0516 &0.0536& 0.0332   & 0.0323 &0.0535& 0.0327 & 0.0436 &0.0575& 0.0325 \\
    & 50 & 0.0480 &0.0532& 0.0393    & 0.0359 &0.0543& 0.0387 & 0.0452 &0.0537& 0.0392 \\
    & 100& 0.0462 &0.0518& 0.0403   & 0.0442 &0.0528& 0.0421 & 0.0445 &0.0518& 0.0421 \\
    & 200& 0.0531 &0.0509& 0.0465   & 0.0507 &0.0524& 0.0473 & 0.0515 &0.0531& 0.0490 \\
  \hline
\end{tabular}
\end{table}
}

{
\begin{table}[p]
\caption{}
\begin{tabular}{ccccc|ccc|ccc}

 \multicolumn{11}{c}{Power of $T_{our}$ and $T_{sk}$ in Case 1} \\
    \hline
\multirow{2}{*}{$p$} &\multirow{2}{*}{$n^*$}&\multicolumn{3}{c}{$z_{ijk}\stackrel{iid}{\sim}N(0,1)$}&\multicolumn{3}{c}{$z_{ijk}\stackrel{iid}{\sim}(\chi^2_2-2)/2$} &\multicolumn{3}{c}{$z_{ijk}\stackrel{iid}{\sim}(\chi^2_8-8)/4$}\\ \cline{3-11}
&&$T_{our}^{umvue}$&$T_{our}^{une}$&$T_{sk}$&$T_{our}^{umvue}$&$T_{our}^{une}$&$T_{sk}$&$T_{our}^{umvue}$&$T_{our}^{une}$&$T_{sk}$\\\hline
 20 & 20 & 0.0941 &0.0985& 0.1079 & 0.0615 &0.0931& 0.1035 & 0.0828 &0.0972& 0.1031 \\
    & 50 & 0.1673 &0.1655& 0.2188 & 0.1380 &0.1630& 0.2295 & 0.1538 &0.1579& 0.2191 \\
    & 100& 0.3111 &0.3205& 0.4649 & 0.3009 &0.3180& 0.4775 & 0.3139 &0.3206& 0.4678 \\
    & 200& 0.6710 &0.9888& 0.8453 & 0.6571 &0.9902& 0.8497 & 0.6664 &0.9888& 0.8529 \\
 50 & 20 & 0.1039 &0.1158& 0.1380 & 0.0748 &0.1130& 0.1494 & 0.0913 &0.1147& 0.1434 \\
     & 50 & 0.2380 &0.2431& 0.3708 & 0.2102 &0.2353& 0.3957 & 0.2244 &0.2372& 0.3808 \\
    & 100& 0.5390 &0.5459& 0.7787 & 0.5218 &0.5419& 0.7879 & 0.5324 &0.5446& 0.7797 \\
    & 200& 0.9342 &1.0000& 0.9918 & 0.9344 &1.0000& 0.9932 & 0.9388 &1.0000& 0.9931 \\
100 & 20 & 0.1280 &0.1404& 0.1843 & 0.0954 &0.1357& 0.2076 & 0.1154 &0.1379& 0.1918 \\
      & 50 & 0.3493 &0.3536& 0.5712 & 0.3183 &0.3577& 0.6043 & 0.3433 &0.3535& 0.5871 \\
     & 100& 0.7787 &0.7707& 0.9574 & 0.7621 &0.7799& 0.9609 & 0.7774 &0.7781& 0.9581 \\
     & 200& 0.9966 &1.0000& 0.9999 & 0.9959 &1.0000& 1.0000 & 0.9965 &1.0000& 0.9999 \\
500 & 20 & 0.2908 &0.3120& 0.4691 & 0.2206 &0.2993& 0.5203 & 0.2723 &0.2989& 0.4775 \\
     & 50 & 0.8576 &0.8690& 0.9884 & 0.8413 &0.8765& 0.9937 & 0.8643 &0.8683& 0.9917 \\
     & 100& 1.0000 &0.9998& 1.0000 & 0.9997 &0.9998& 1.0000 & 1.0000 &0.9999& 1.0000 \\
     & 200& 1.0000 &1.0000& 1.0000 & 1.0000 &1.0000& 1.0000 & 1.0000 &1.0000& 1.0000 \\
800 & 20 & 0.3921 &0.4103& 0.6118 & 0.3028 &0.4119& 0.6763 & 0.3706 &0.4001& 0.6312 \\
     & 50 & 0.9670 &0.9662& 0.9997 & 0.9547 &0.9668& 0.9997 & 0.9641 &0.9681& 0.9999 \\
     & 100& 1.0000 &1.0000& 1.0000 & 1.0000 &1.0000& 1.0000 & 1.0000 &1.0000& 1.0000 \\
     & 200& 1.0000 &1.0000& 1.0000 & 1.0000 &1.0000& 1.0000 & 1.0000 &1.0000& 1.0000 \\
  \hline
\end{tabular}
\end{table}
}

{
\begin{table}[p]
\caption{}
\begin{tabular}{ccccc|ccc|ccc}

 \multicolumn{11}{c}{ASL of $T_{our}$ and $T_{sk}$ in Case 2} \\
    \hline
\multirow{2}{*}{$p$} &\multirow{2}{*}{$n^*$}&\multicolumn{3}{c}{$z_{ijk}\stackrel{iid}{\sim}N(0,1)$}&\multicolumn{3}{c}{$z_{ijk}\stackrel{iid}{\sim}(\chi^2_2-2)/2$} &\multicolumn{3}{c}{$z_{ijk}\stackrel{iid}{\sim}(\chi^2_8-8)/4$}\\ \cline{3-11}
&&$T_{our}^{umvue}$&$T_{our}^{une}$&$T_{sk}$&$T_{our}^{umvue}$&$T_{our}^{une}$&$T_{sk}$&$T_{our}^{umvue}$&$T_{our}^{une}$&$T_{sk}$\\\hline
 20 & 20 & 0.0723 &0.0707& 0.0153 & 0.0679 &0.0788& 0.0103 & 0.0739 &0.0705& 0.0120 \\
    & 50 & 0.0688 &0.0673& 0.0087 & 0.0686 &0.0657& 0.0105 & 0.0654 &0.0697& 0.0084 \\
    & 100& 0.0708 &0.0695& 0.0100 & 0.0705 &0.0730& 0.0085 & 0.0664 &0.0670& 0.0096 \\
    & 200& 0.0695 &0.0656& 0.0083 & 0.0622 &0.0710& 0.0061 & 0.0640 &0.0664& 0.0104 \\
 50 & 20 & 0.0755 &0.0747& 0.0098 & 0.0815 &0.0736& 0.0116 & 0.0646 &0.0755& 0.0092 \\
    & 50 & 0.0711 &0.0722& 0.0061 & 0.0676 &0.0725& 0.0069 & 0.0685 &0.0735& 0.0066 \\
    & 100& 0.0683 &0.0694& 0.0074 & 0.0623 &0.0707& 0.0052 & 0.0658 &0.0695& 0.0070 \\
    & 200& 0.0687 &0.0688& 0.0064 & 0.0673 &0.0645& 0.0066 & 0.0645 &0.0709& 0.0059 \\
100 & 20 & 0.0744 &0.0747& 0.0086 & 0.0710 &0.0705& 0.0067 & 0.0745 &0.0722& 0.0074 \\
    & 50 & 0.0687 &0.0720& 0.0067 & 0.0674 &0.0665& 0.0056 & 0.0744 &0.0724& 0.0072 \\
    & 100& 0.0679 &0.0719& 0.0056 & 0.0649 &0.0641& 0.0046 & 0.0752 &0.0690& 0.0067 \\
    & 200& 0.0699 &0.0705& 0.0057 & 0.0695 &0.0684& 0.0052 & 0.0673 &0.0674& 0.0056 \\
500 & 20 & 0.0736 &0.0710& 0.0044 & 0.0755 &0.0722& 0.0048 & 0.0690 &0.0754& 0.0043 \\
    & 50 & 0.0708 &0.0714& 0.0043 & 0.0703 &0.0700& 0.0043 & 0.0656 &0.0681& 0.0032 \\
    & 100& 0.0703 &0.0705& 0.0035 & 0.0713 &0.0613& 0.0019 & 0.0660 &0.0728& 0.0023 \\
    & 200& 0.0694 &0.0692& 0.0032 & 0.0686 &0.0680& 0.0030 & 0.0694 &0.0709& 0.0019 \\
800 & 20 & 0.0730 &0.0753& 0.0041 & 0.0797 &0.0797& 0.0037 & 0.0736 &0.0767& 0.0032 \\
    & 50 & 0.0680 &0.0697& 0.0025 & 0.0686 &0.0738& 0.0026 & 0.0754 &0.0733& 0.0026 \\
    & 100& 0.0686 &0.0667& 0.0022 & 0.0667 &0.0651& 0.0017 & 0.0680 &0.0683& 0.0017 \\
    & 200& 0.0705 &0.0686& 0.0015 & 0.0688 &0.0678& 0.0028 & 0.0691 &0.0620& 0.0016 \\
  \hline
\end{tabular}
\end{table}
}

{
\begin{table}[p]
\caption{}
\begin{tabular}{ccccc|ccc|ccc}

 \multicolumn{11}{c}{Power of $T_{our}$ and $T_{sk}$ in Case 2} \\
    \hline
\multirow{2}{*}{$p$} &\multirow{2}{*}{$n^*$}&\multicolumn{3}{c}{$z_{ijk}\stackrel{iid}{\sim}N(0,1)$}&\multicolumn{3}{c}{$z_{ijk}\stackrel{iid}{\sim}(\chi^2_2-2)/2$} &\multicolumn{3}{c}{$z_{ijk}\stackrel{iid}{\sim}(\chi^2_8-8)/4$}\\ \cline{3-11}
&&$T_{our}^{umvue}$&$T_{our}^{une}$&$T_{sk}$&$T_{our}^{umvue}$&$T_{our}^{une}$&$T_{sk}$&$T_{our}^{umvue}$&$T_{our}^{une}$&$T_{sk}$\\\hline
 20 & 20 & 0.1952 &0.0998& 0.0573 & 0.1855 &0.1998& 0.0529 & 0.1859 &0.1962&0.0557 \\
     & 50 & 0.3860 &0.1420& 0.1377 & 0.3814 &0.3973& 0.1345 & 0.3927 &0.3881& 0.1406 \\
    & 100& 0.6769 &0.2225& 0.3454 & 0.6783 &0.6782& 0.3486 & 0.6864 &0.6799& 0.3457 \\
    & 200& 0.9464 &0.7934& 0.7498 & 0.9491 &0.9999& 0.7450 & 0.9516 &0.9998& 0.7475 \\
 50 & 20 & 0.2096 &0.1111& 0.0486 & 0.2031 &0.2215& 0.0500 & 0.2016 &0.2087& 0.0502 \\
     & 50 & 0.4201 &0.1583& 0.1327 & 0.4246 &0.4257& 0.1326 & 0.4350 &0.4329& 0.1385 \\
    & 100& 0.7369 &0.2442& 0.3462 & 0.7425 &0.7423& 0.3588 & 0.7506 &0.7471& 0.3583 \\
    & 200& 0.9779 &0.8457& 0.7835 & 0.9796 &1.0000& 0.7763 & 0.9816 &1.0000& 0.7899 \\
100 & 20 & 0.2148 &0.1050& 0.0471 & 0.2112 &0.2230& 0.0471 & 0.2163 &0.2188& 0.0458 \\
    & 50 & 0.4575 &0.1529& 0.1286 & 0.4564 &0.4583& 0.1305 & 0.4631 &0.4571& 0.1292 \\
    & 100& 0.7860 &0.2482& 0.3601 & 0.7838 &0.7839& 0.3549 & 0.7797 &0.7863& 0.3482  \\
    & 200& 0.9899 &0.8846& 0.7939 & 0.9897 &1.0000& 0.7994 & 0.9912 &1.0000& 0.7949 \\
500 & 20 & 0.2465 &0.1134& 0.0342 & 0.2430 &0.2478& 0.0331 & 0.2450 &0.2460& 0.0350 \\
    & 50 & 0.5226 &0.1787& 0.1127 & 0.5220 &0.5162& 0.1119 & 0.5284 &0.5198& 0.1161 \\
    & 100& 0.8594 &0.2807& 0.3373 & 0.8624 &0.8637& 0.3371 & 0.8595 &0.8619& 0.3381 \\
    & 200& 0.9990 &0.9400& 0.8238 & 0.9995 &1.0000& 0.8171 & 0.9993 &1.0000& 0.8212 \\
800 & 20 & 0.2474 &0.1237& 0.0290 & 0.2586 &0.2572& 0.0320 & 0.2579 &0.2597& 0.0317 \\
    & 50 & 0.5397 &0.1867& 0.1063 & 0.5550 &0.5482& 0.1061 & 0.5466 &0.5352& 0.1033 \\
    & 100& 0.8250 &0.2984& 0.3359 & 0.8850 &0.8893& 0.3459 & 0.8855 &0.8865& 0.3370 \\
    & 200& 0.9998 &0.9548& 0.8152 & 0.9997 &1.0000& 0.8252 & 0.9996 &1.0000& 0.8247 \\
  \hline
\end{tabular}
\end{table}
}

\section{Technical details}
In this section we give the proof of Theorem \ref{th1}. We restricted our
attention to the case in which $k = 3$ for simplicity  and the proof for the case of $k>3$ is the same.
Here we use the same method as in \cite{ChenQ10T}, hence some of the derivations are omitted. The main difference is that we need to verify the asymptotic normality of $T_n^{(3)}$.  Because it does not follow by any means that the random variable $\a_n+\b_n$ will converge in distribution to $\a+\b$, if $\a_n\stackrel{d}{\to}\a$ and $\b_n\stackrel{d}{\to}\b$.

Denote $T_n^{(3)}=T_{n1}^{(3)}+T_{n2}^{(3)}$, where
\begin{align*}
  T^{(3)}_{n1}=2\sum_{k=1}^3\sum_{i\not=j}^{n_k}\frac{(X_{ki}-\mu_k)'(X_{kj}-\mu_k)}{n_k(n_k-1)}
-2\sum_{k<l}^3\sum_{i=1}^{n_k}\sum_{j=1}^{n_l}\frac{(X_{ki}-\mu_k)'(X_{lj}-\mu_l)}{n_kn_l}
\end{align*}
and
\begin{align*}
   T^{(3)}_{n2}=2\sum_{k,l=1}^3\sum_{i=1}^{n_k}\frac{(X_{ki}-\mu_k)'(\mu_k-\mu_l)}{n_k}+\sum_{k<l}^3\parallel\mu_k-\mu_l\parallel^2
\end{align*}
We can verify that $E(T^{(3)}_{n1})=0$, $E(T^{(3)}_{n2})=\sum_{k<l}^3\parallel\mu_k-\mu_l\parallel^2$ and
\begin{align*}
Var(T^{(3)}_{n2})=&4\sum_{k<l}^3(\mu_k-\mu_l)'(n_l^{-1}\Sigma_l+n_k^{-1}\Sigma_k)(\mu_k-\mu_l)\\&+4\sum_{i\neq j\neq k }^3(\mu_i-\mu_j)'n_i^{-1}\Si_i(\mu_i-\mu_k).
\end{align*}
From condition (e), that is,
\begin{align*}
  Var\(\frac{ T^{(3)}_{n2}-\sum_{k<l}^3\parallel\mu_k-\mu_l\parallel^2}{\si_n^{(3)}}\)=o(1),
\end{align*}
we get
\begin{align*}
  \frac{ T^{(3)}_{n}-\sum_{k<l}^3\parallel\mu_k-\mu_l\parallel^2}{\sqrt {Var(T_n^{(3)})}}=\frac{ T^{(3)}_{n1}}{\si_n^{(3)}}+o_p(1).
\end{align*}

Next we will prove  the asymptotic normality of $T_{n1}^{(3)}$.  Without loss of generality we assume that $\mu_1=\mu_2=\mu_3=0$. Let $Y_i=X_{1i}(i=1,\ldots,n_1)$, $Y_{j+n_1}=X_{2j}(j=1,\ldots,n_2)$, $Y_{j+n_1+n_2}=X_{3j}(j=1,\ldots,n_3)$. For  $i\not=j$,
\begin{displaymath}
\phi_{ij} = \left\{ \begin{array}{ll}
2n_1^{-1}(n_1-1)^{-1}Y'_i Y_j, & \textrm{if $i,j\in\{1,2,\ldots,n_1\}$};\\
-n_1^{-1}n_2^{-1}Y'_i Y_j, & \textrm{if $i\in\{1,2,\ldots,n_1\}$ and $j\in\{n_1+1,\ldots,n_1+n_2\}$};\\
2n_2^{-1}(n_2-1)^{-1}Y'_i Y_j, & \textrm{if $i,j\in\{n_1+1,\ldots,n_1+n_2\}$};\\
-n_2^{-1}n_3^{-1}Y'_i Y_j, & \textrm{if $i\in\{n_1+1,\ldots,n_1+n_2\}$ }\\
&\qquad\quad\textrm{and $j\in\{n_1+n_2+1,\ldots,n_1+n_2+n_3\}$};\\
2n_3^{-1}(n_3-1)^{-1}Y'_i Y_j, & \textrm{if $i,j\in\{n_1+n_2+1,\ldots,n_1+n_2+n_3\}$};\\
-n_3^{-1}n_1^{-1}Y'_i Y_j, & \textrm{if $i\in\{1,2,\ldots,n_1\}$ and $j\in\{n_1+n_2+1,\ldots,n_1+n_2+n_3\}$}.\\
\end{array} \right.
\end{displaymath}
For $j=2,3,\ldots,n_1+n_2+n_3$, denote $V_{nj}=\sum_{i=1}^{j-1}\phi_{ij}$, $S_{nm}=\sum_{j=2}^{m}V_{nj}$  and
$\mathcal{F}_{nm}=\sigma\{Y_1,Y_2,\ldots,Y_m\}$ which is the $\si $ algebra generated by $\{Y_1,Y_2,\ldots,Y_m\}$. Then we have
$$T_n^{(3)}=2\sum_{j=2}^{n_1+n_2+n_3} V_{nj}.$$
 It is easy to verify that $\{S_{nm},\mathcal{F}_{nm}\}_{m=1}^n$  forms a sequence of zero mean and 
square integrable martingale.
Then the asymptotic normality of $T_n^{(3)}$ can  be proved by employing Corollary 3.1 in
 \cite{HallH80M} with routine verification of the following:
\begin{align}\lb{ma1}
  \frac{\sum_{j=2}^{n_1+n_2+n_3}E[V_{nj}^2|\mathcal{F}_{n,j-1}]}{(\si_n^{(3)})^2}\stackrel{P}{\to}\frac{1}{4}.
\end{align}
and
\begin{align}\lb{ma2}
\sum_{j=2}^{n_1+n_2+n_3}(\si_n^{(3)})^{-2}E[V^2_{nj}I(|V_{nj}|>\epsilon\si_n^{(3)})
|\mathcal{F}_{n,j-1}]\stackrel{P}{\to}0.
\end{align}
Thus next  we prove \eqref{ma1} and \eqref{ma2} respectively.
\subsection{Proof of \eqref{ma1}.}
Verify that $$E(V_{nj}^2|\mathcal{F}_{n,j-1})=E\Big[\Big(\sum_{i=1}^{j-1}\phi_{ij}\Big)^2\Big|\mathcal{F}_{n,j-1}\Big]
=E\Big(\sum_{i_1,i_2=1}^{j-1}\phi_{i_1j}\phi_{i_2j}\Big|\mathcal{F}_{n,j-1}\Big)$$
$$=\sum_{i_1,i_2=1}^{j-1}c_{i_1j}c_{i_2j}Y'_{i_1}E(Y_jY'_j|\mathcal{F}_{n,j-1})Y_{i_2}=\sum_{i_1,i_2=1}^{j-1}c_{i_1j}c_{i_2j}
Y'_{i_1}E(Y_jY'_j)Y_{i_2}$$
$$=\sum_{i_1,i_2=1}^{j-1}c_{i_1j}c_{i_2j}Y'_{i_1}{\tilde{\Sigma}_j}Y_{i_2},
$$
where $c_{ij}$ is the coefficient of $\phi_{ij}$ and if $j\in[1,n_1]$, $\tilde{\Sigma}_j=\Sigma_1$;  if $j\in[n_1+1,n_1+n_2]$,
$\tilde{\Sigma}_j=\Sigma_2$; if $j\in[n_1+n_2+1,n_1+n_2+n_3]$, $\tilde{\Sigma}_j=\Sigma_3$.

Denote $$
\eta_n=\sum_{j=2}^{n_1+n_2+n_3}E(V_{nj}^2|\mathcal{F}_{n,j-1}).$$
Then we have
\begin{align*}
E(\eta_n)=&\frac{2tr(\Sigma_1^2)}{n_1(n_1-1)}+\frac{2tr(\Sigma_2^2)}{n_2(n_2-1)}+\frac{2tr(\Sigma_3^2)}{n_3(n_3-1)}\\
&+
\frac{tr(\Sigma_1\Sigma_2)}{n_1n_2}+\frac{tr(\Sigma_1\Sigma_3)}{n_1n_3}+\frac{tr(\Sigma_3\Sigma_2)}{n_3n_2}=\frac{1}{4}(\sigma_{n}^{(3)})^2.
\end{align*}

Now consider
\begin{align}\lb{eqa=}
E(\eta_n^2)=E\Big[\sum_{j=2}^{n_1+n_2+n_3}\sum_{i_1,i_2=1}^{j-1} c_{i_1j}c_{i_2j}Y'_{i_1}{\tilde{\Sigma}_j}Y_{i_2}\Big]^2=2E(A)+E(B),
\end{align}
where
\begin{align*}
A=\sum_{2\leq j_1<j_2}^{n_1+n_2+n_3}\sum_{i_1,i_2=1}^{{j_1}-1}\sum_{i_3,i_4=1}^{{j_2}-1}
c_{i_1j_1}c_{i_2j_1}c_{i_3j_2}c_{i_4j_2}Y'_{i_1}{\tilde{\Sigma}_{j_1}}Y_{i_2}Y'_{i_3}{\tilde{\Sigma}_{j_2}}Y_{i_4}\\
\end{align*}and
\begin{align*}
B=\sum_{j=2}^{n_1+n_2+n_3}\sum_{i_1,i_2=1}^{j-1}\sum_{i_3,i_4=1}^{j-1}
c_{i_1j}c_{i_2j}c_{i_3j}c_{i_4j}Y'_{i_1}{\tilde{\Sigma}_j}Y_{i_2}Y'_{i_3}{\tilde{\Sigma}_j}Y_{i_4}.
\end{align*}
The term $B$ can be further partitioned as $B = B_1 + B_2+B_3$, where
\begin{align*}
  E(B_1)=&E\(\sum_{j=2}^{n_1}\sum_{i_1,i_2=1}^{j-1}\sum_{i_3,i_4=1}^{j-1}
c_{i_1j}c_{i_2j}c_{i_3j}c_{i_4j}Y'_{i_1}{{\Sigma}_1}Y_{i_2}Y'_{i_3}{{\Sigma}_1}Y_{i_4}\)\\
  E(B_2)=&E\(\sum_{j=n_1+1}^{n_1+n_2}\sum_{i_1,i_2=1}^{j-1}\sum_{i_3,i_4=1}^{j-1}
c_{i_1j}c_{i_2j}c_{i_3j}c_{i_4j}Y'_{i_1}{{\Sigma}_2}Y_{i_2}Y'_{i_3}{{\Sigma}_2}Y_{i_4}\)\\
 E(B_3)=&E\(\sum_{j=n_1+n_2+1}^{n_1+n_2+n_3}\sum_{i_1,i_2=1}^{j-1}\sum_{i_3,i_4=1}^{j-1}
c_{i_1j}c_{i_2j}c_{i_3j}c_{i_4j}Y'_{i_1}{{\Sigma}_3}Y_{i_2}Y'_{i_3}{{\Sigma}_3}Y_{i_4}\).
\end{align*}
We only compute $E(B_3 )$ here as $E(B_1 )$ and $E(B_2 )$ can be computed following the same procedure.
As  $\mu_1 = \mu_2 =\mu_3= 0$, we only need to consider $i_1$, $i_2$, $i_3$ and $i_4$ in these three cases:  (a)
$(i_1 =  i_2 )\neq  (i_3 =  i_4 )$;  (b) $(i_1 =  i_3 )\neq  (i_2 =  i_4 )$ or $(i_1 =  i_4 )\neq  (i_2 =  i_3 )$;  (c)
$i_1 =  i_2 =i_3 =  i_4 $.  Thus we obtain that
\begin{align*}
  E(B_3)=E(B_{31})+E(B_{32})+E(B_{33}),
\end{align*}
where
\begin{align*}
 E(B_{31})=&O({n^{-8}})E\(\sum_{j={n_1}+n_2+1}^{{n_1}+n_2+n_3}\sum_{i_1\neq i_2}^{j-1}
Y'_{i_1}{\Si_3}Y_{i_1}Y'_{i_2}{\Si_3}Y_{i_2}\)
=O({n^{-5}})\sum_{i,j=1}^3tr\Si_3\Si_itr\Si_3\Si_j\\
 E(B_{32})=&O({n^{-8}})E\(\sum_{j={n_1}+n_2+1}^{n_1+n_2+n_3}\sum_{i_1\neq i_2}^{j-1}
Y'_{i_1}{\Si_3}Y_{i_2}Y'_{i_2}{\Si_3}Y_{i_1}\)
=O({n^{-5}})\sum_{i,j=1}^3tr\Si_i\Si_3\Si_j\Si_3
\end{align*}
and
\begin{align*}
 E(B_{33})=&O({n^{-8}})E\(\sum_{j=n_1+n_2+1}^{n_1+n_2+n_3}\sum_{i=1}^{j-1}
Y'_{i}{\Si_3}Y_{i}Y'_{i}{\Si_3}Y_{i}\)\\
=&O({n^{-6}})\(\sum_{i=1}^{3}\(E(Z_{i1}'\Ga_i'{\Si_3}\Ga_iZ_{i1}-tr\Si_3\Si_i)^2+tr^2\Si_3\Si_i\)\)\\
=&O({n^{-6}})\sum_{i=1}^{3}\(tr(\Si_3\Si_i)^2+tr^2\Si_3\Si_i\).
\end{align*}
Thus  we obtain that
\begin{align*}
  E(B_3)=o((\sigma_{n}^{(3)})^4).
\end{align*}
As we can similarly get $E(B_1)=o((\sigma_{n}^{(3)})^4)$ and $E(B_2)=o((\sigma_{n}^{(3)})^4)$,  we conclude that
\begin{align}\lb{eb}
    E(B)=o((\sigma_{n}^{(3)})^4).
\end{align}

Using the same method of deriving \eqref{eb}, we have $$2E(A)=\frac1{16}(\sigma_{n}^{(3)})^4(1+o(1)),$$
which together with \eqref{eqa=} and \eqref{eb} implies $$
E(\eta^2_n)=\frac1{16}(\sigma_{n}^{(3)})^4+o((\sigma_{n}^{(3)})^4).$$
Then we have $$Var(\eta_n)=E(\eta_n^2)-E^2(\eta_n)=o((\sigma_{n}^{(3)})^4).$$
Therefore we obtain
$$(\sigma_{n}^{(3)})^{-2}E\Big\{\sum_{j=1}^{n_1+n_2+n_3}E(V^2_{nj}|\mathcal{F}_{n,j-1})\Big\}
=(\sigma_{n}^{(3)})^{-2}E(\eta_n)=\frac{1}{4}$$
and $$(\sigma_{n}^{(3)})^{-4}Var\Big\{\sum_{j=1}^{n_1+n_2+n_3}E(V^2_{nj}|\mathcal{F}_{n,j-1})\Big\}
=(\sigma_{n}^{(3)})^{-4}Var(\eta_n)=o(1),$$
which complete the proof of  \eqref{ma1}.
\subsection{Proof of \eqref{ma2}.}
As $$\sum_{j=2}^{n_1+n_2+n_3}(\sigma_{n}^{(3)})^{-2}E\{V^2_{nj}I(|V_{nj}|>\epsilon{\sigma_{n}^{(3)})
|\mathcal{F}}_{n,j-1}\}\leq(\sigma_{n}^{(3)})^{-4}\epsilon^{-2}\sum_{j=1}^{n_1+n_2+n_3}E(V^4_{nj}|{\mathcal{F}}_{n,j-1}),$$
we just need to show that
$$E\(\sum_{j=2}^{n_1+n_2+n_3}E(V^4_{nj}|{\mathcal{F}}_{n,j-1})\)=o((\sigma^{(3)}_{n})^4).$$
{Note} that
\begin{align*}
&E\Big\{\sum_{j=2}^{n_1+n_2+n_3}E(V^4_{nj}|{\mathcal{F}}_{n,j-1})\Big\}=
\sum_{j=2}^{n_1+n_2+n_3}E(V^4_{nj})=O(n^{-8})\sum_{j=2}^{n_1+n_2+n_3}E\Big(
\sum_{i=1}^{j-1}Y_{i}'Y_j\Big)^4\\
=&O(n^{-8})\sum_{j=2}^{n_1+n_2+n_3}\sum_{s\not=t}^{j-1}E((Y'_jY_s)^2(Y'_tY_j)^2) +O(n^{-8})\sum_{j=2}^{n_1+n_2+n_3}\sum_{s=1}^{j-1}E(Y'_sY_j)^4\\
=&O(n^{-8})\sum_{j=2}^{n_1+n_2+n_3}\sum_{s\not=t}^{j-1}E(Y'_j\tilde{\Sigma}_sY_jY'_j\tilde{\Sigma}_tY_j) +O(n^{-8})\sum_{j=2}^{n_1+n_2+n_3}\sum_{s=1}^{j-1}E(Y'_sY_j)^4.
\end{align*}
The first term of last equation has the order $o((\sigma^{(3)}_{n})^4)$ which can be proved by the {same} procedure in last subsection. It remains to consider the second term. As proved in \cite{ChenQ10T}, we have
\begin{align*}
 \sum_{j=2}^{n_1+n_2}\sum_{s=1}^{j-1}E(Y'_sY_j)^4
  =O({n^2})\(\sum_{i=1}^2(tr^2(\Si_i^2)+tr(\Si_i^4))+tr^2(\Si_1\Si_2)+tr(\Si_1\Si_2)^2\),
\end{align*}
and
\begin{align*}
  &\sum_{j=n_1+n_2+1}^{n_1+n_2+n_3}\sum_{s=1}^{j-1}E(Y'_sY_j)^4 =\sum_{j=n_1+n_2+1}^{n_1+n_2+n_3}\sum_{s=1}^{n_1}E(Y'_sY_j)^4 \\ &+\sum_{j=n_1+n_2+1}^{n_1+n_2+n_3}\sum_{s=n_1+1}^{n_1+n_2}E(Y'_sY_j)^4 +\sum_{j=n_1+n_2+1}^{n_1+n_2+n_3}\sum_{s=n_1+n_2+1}^{j-1}E(Y'_sY_j)^4\\
  =&{O(n^2)\(tr^2(\Si_3^2)+tr(\Si_3^4)+\sum_{i=1}^2tr^2(\Si_i\Si_3)+\sum_{i=1}^2tr(\Si_i\Si_3)^2\)}.
  \end{align*}
Thus
we conclude that
\begin{align*}
 {O(n^{-8})}\sum_{j=2}^{n_1+n_2+n_3}\sum_{s=1}^{j-1}E(Y'_sY_j)^4= {O(n^{-8})}\sum_{j=2}^{{n_1+n_2}}\sum_{s=1}^{j-1}E(Y'_sY_j)^4+ {O(n^{-8})}\sum_{j={n_1+n_2+1}}^{n_1+n_2+n_3}\sum_{s=1}^{j-1}E(Y'_sY_j)^4\\
{=O(n^{-6})\(\sum_{i=1}^3(tr^2(\Si_i^2)+tr(\Si_i^4))+\sum_{i<j}^3tr^2(\Si_i\Si_j)+\sum_{i<j}^3tr(\Si_i\Si_j)^2\)}=o((\sigma^{(3)}_{n})^4).
\end{align*}
Then the proof of \eqref{ma2} is  complete.

\section*{Acknowledgements}
The authors would like to thank the referees for many constructive comments. 


\end{document}